\documentclass{amsart}
\usepackage{amssymb}
\newtheorem{thm}{Theorem}[section]
\newtheorem{lma}{Lemma}[section]

\newtheorem*{definition}{Definition}
\theoremstyle{remark}
\newtheorem*{acknowledgement}{Acknowledgement}

\newcommand{\md}{\medskip}
\newcommand{\al}{\alpha}
\numberwithin{equation}{section}

\begin{document}
\title[periodic solutions of Duffing
equations]{Stability and exact multiplicity of\\periodic solutions of Duffing
equations\\with cubic nonlinearities}
\author{Hongbin Chen}
\address{Department of Mathematics, Xi'an Jiaotong
University, Xi'an P.R.~China}
\email{hbchen@mail.xjtu.edu.cn}
\author{Yi Li}
\address{Department of Mathematics, The University of Iowa, Iowa City, IA 52242, USA
and
Department of Mathematics, Hunan Normal University, Changsha, Hunan, China}

\begin{abstract}
We study the stability and exact multiplicity of periodic solutions
of the Duffing equation with cubic nonlinearities,
$$
x''+cx'+ax-x^{3}=h(t), \eqno {(*)}
$$
 where $a$ and $c>0$ are positive constants and $h(t)$ is a positive $T$-periodic function.
 We obtain sharp bounds for $h$ such that $(*)$ has exactly three ordered $T$-periodic
 solutions. Moreover, when $h$ is within these bounds, one of the three solutions is negative, while the
 other two are positive. The middle solution is asymptotically stable, and the remaining two are unstable.
\end{abstract}
\maketitle

\section{Introduction}

Consider the  Duffing equation
 \begin{equation} \label {aa}x''+cx'+ax-x^{3}=h(t),\quad x(0)=x(T),\quad x'(0)=x'(T),\end{equation}
where $h$ is a positive
$T$-periodic function and $a$ and $c$ are
constants with $0<a<(\frac{\pi}{T})^{2}+\frac{c^{2}}{4}$ and $c>0$.
 The reason that we assume $0<a<(\frac{\pi}{T})^{2}+\frac{c^{2}}{4}$ is that for $0>a$,
(\ref{aa})
 has a unique $T$-periodic solution: there is no bifurcation at all. Therefore the only
 interesting case is $a>0$. We assume that
 \begin{equation}\label{ab} a<(\frac{\pi}{T})^{2}+\frac{c^{2}}{4}\end{equation}
  Namely, the system is
 sufficient damped probable due to the drawback of our proof, we
 still don't know whether the condition(\ref{ab}) is essential or
 not.
 The existence and multiplicity of
periodic solutions of (\ref{aa}) or more general types of
nonlinear second-order differential equation have been
investigated extensively by many authors.  However, relatively few
studies have been written about the exact multiplicity
 of (\ref{aa}). In \cite{ChenL}, we have studied the small-perturbation problem, and
 established that (\ref{aa}) has exactly three $T$-periodic solutions provided
 that $h(t)$ is sufficiently small.
  R.~Ortega has considered the following parametrized
Duffing equation\cite%
{Ort}:
\begin{equation}
 \label {ac}x''+cx'+g(t,x)= s+h(t),
\end{equation}
where $c>0$ and $
g'_{x}(t,x)$ is strictly increasing in the second variable.
Under the additional assumptions
\[
g'_{x}(t,x)\ll \frac{{\pi}^{2}}{T^{2}}+\frac{c^2}{4}
\]
and
\[
\lim_{x\rightarrow \pm\infty}g(t,x)=+\infty,
\]
he has shown the following Ambrosetti-Prodi-type theorem:\md
\begin{quotation}
There is an $s_{0}$ such that (\ref{ac}) has no $T$-periodic
solution for $s<s_{0}$, (\ref{ac}) has a unique $T$-periodic
solution which is unstable for $s=s_{0}$, and (\ref{ac}) has exactly
two periodic solutions for $s>s_{0}$, one of them is asymptotically
stable and another  is unstable.\md
\end{quotation}
For related results on the existence of two periodic solutions, we refer to
\cite{FMN,Tineo,Ga}. Similar results concerning the first-order equation
were obtained by J.~Mawhin \cite{Maw}, and more recently by the authors \cite%
{Chen} based on singularity theory and A.~Tineo \cite{Tineoa}. For the
multiplicity results concerning the forced pendulum equation, one can refer to
\cite{FMa,Katr,Tara}. However, as far as the authors know, there is no such
precise result on existence of exactly three periodic solutions for nonlinear
Duffing equations. Here, we obtain multiplicity for periodic solutions by
means of a topological degree argument combined with a newly developed maximum
principle given in \cite{Torres,Omar,Njo}, and stability of periodic
solutions follows by computing the local index given by R.~Ortega in \cite%
{Ort}. The more recent results concerning the stability and the sharpness of
rate of decay of periodic solutions can be found in \cite%
{Al,Alort,ChenLi,LM,LaMc}.

\begin{thm}
\label{Theorema} Let $h_{0}:=\sqrt{4a^{3}/27}$. Then
\begin{enumerate}
\item\label{thm:1} \textup{(\ref {aa})} has a unique  $T$-periodic solution
 which is negative and unstable if $h(t)>h_{0} \;\forall\, t \in \mathbb{R}$;
\item \label{thm:2}  \textup{(\ref{aa})} has  exactly three  ordered
$T$-periodic solutions if $0<h(t)<h_{0}$.
 \item\label{thm:3}Moreover, in case \textup{(\ref{thm:2})}, the minimal
solution is negative and the other two are positive; also, the middle solution is
asymptotically stable and the remaining two are unstable.

\end{enumerate}
\end{thm}

  The following notations are used.
\begin{enumerate}
\item $L_{T }^{p}$\qquad $T $-periodic function $u\in L^{p}[0,T]$ with $
\Vert u\Vert _{p}$ for $1\leq p\leq \infty$;
\item $C_{T }^{k}$\qquad $T $-periodic function $u\in C^{k}[0,T
]$, $k\geq 0$, with $C^{k}$-norm;
\item $\al(t)\gg\beta(t)$, if $
\al(t)\geq\beta(t)$ and $ \al(t)>\beta(t)$ on some positive-%
measure subset.
\end{enumerate}

\section{Topological index and linear periodic problems}

In this section we shall recall some basic results about linear
 periodic boundary-value problems  that will be needed in the sequel.

 Consider the periodic boundary-value
problem
\begin{equation}
\label{tla}  \left\{
\begin{aligned}
& x'=F(t,x),\\
& x(0)=x(T),
\end{aligned}
\right.
\end{equation}
where $F\colon  [0,T ]\times {\mathbb{R}^{n}}\rightarrow
 \mathbb{R}^{n}$ is a continuous function that is $ T $-periodic in
$t$ for $n=2$.

We denote by $x(t,x_{0})$
the initial-value solution of (\ref{tla}).

\begin{definition}A
$T$-periodic solution $x$ of \textup{(\ref{tla})} will be called a
\emph{nondegenerate} $T $-periodic solution if the linearized equation
\begin{equation}
 y^{\prime }=F_{x}(t,x)y \label{eqtb}
\end{equation}
 does not admit a
nontrivial $T $-periodic solution.
\end{definition}

Let $M(t)$ be the fundamental matrix of (\ref{eqtb}), and $\mu
_{1}$ and $\mu _{2}$ the eigenvalues of the matrix $M(T )$. Then
$x(t,x_{0})$ is asymptotically stable if and only if $\vert \mu
_{i}\vert <1$, $i=1,2$; otherwise, if one of them has modulus
greater than one, $x(t,x_{0})$ is unstable.

 Consider the homogeneous
periodic equation
\begin{equation}
\label{eqtc} L_{\al }x=x^{\prime \prime
}+cx^{\prime }+\al (t)x=0
\end{equation}
 where $c$ is a constant and $\al
(t)\in L_{T }$.

The following simple lemma is given by the authors in \cite%
{ChenL} which is needed  in proving our main results.

\begin{lma}
\label{tlma}Suppose that $\al (t)$, $\al _{1}(t)$ and $\al
_{2}(t)\in L_{T }$ such that
\begin{equation}
\label{eqtx} \al _{1}(t),\;\al
_{2}(t)\textup{ and }\al (t)\ll
\left(\frac{2\pi}{T}\right)^{2}+\frac{c^{2}}{4}.
\end{equation}
Then
\begin{enumerate}
\item \label{l(1)}the possible $T $-periodic solution $x$ of
equation \textup{(\ref{eqtc})} is either trivial or different from
zero for each $t\in \mathbb{R}$; \item \label{l(2)} $L_{\al
_{i}}x=0$ \textup{(}$i=1,2$\textup{)} cannot admit nontrivial
$T$-periodic solutions simultaneously if $\al _{1}(t)\ll \al
_{2}(t)$; \item \label{l(3)}$L_{\al }x=0$ has no nontrivial
$T$-periodic solution, and $\operatorname{ind} L_{\al}:=\deg(L_{\al
})=1$ \textup{(}resp., $=-1$\textup{)}, if $\al (t)\gg 0$
\textup{(}resp., $\al(t)\ll 0$\textup{)}.
\end{enumerate}
\end{lma}

The following connection between stability and topological index
is due to R.~Ortega \cite{Ort}.

\begin{lma}
\label{lmata} Assume that $x$ is an isolated $T$-periodic solution
of \textup{(\ref{ab})} such that the condition
\[
g_{x}(t,x)\leq \left(\frac{\pi}{T}\right)^{2}+\frac{c^{2}}{4}
\]
holds, for
$t\in \mathbb{R}$, $c>0$. Then $x$ is asymptotically stable \textup{(}resp.,
unstable\/\textup{)} if and only if $\operatorname{ind}(x)=1$ \textup{(}resp.,
$\operatorname{ind}(x)=-1$\textup{)}.
\end{lma}

Consider the differential equation
\begin{equation}
\label{eqte} x^{\prime \prime }+cx^{\prime }+\al (t)x=h(t),
\end{equation}
 where $c$ is
a constant and $\al (t)$, $h(t)\in L_{T }$.

The following  maximum principle is given by P.J. Torres and M.R.
Zhang (Theorem 2.3) in \cite{Torres}. Here we state the principle in
a somewhat different form.
\begin{lma}
\label{lae} Let $h(t)\gg0$ and $\al(t)$ satisfy
\[
\al(t)\leq \left(\frac{\pi}{T}\right)^{2}+\frac{c^{2}}{4}.
\]
If $x(t)$ is a
$T$-periodic solution of \textup{(\ref{eqte})}, then the following
statements hold:
\begin{enumerate}
\item \label{la} Either $x(t)>0$ or $x(t)<0$ for all $t \in
\mathbb{R}$; \item \label{lb} $x(t)>0 \;\forall\, t \in \mathbb{R}$
if $\al(t)>0\; \forall\, t \in \mathbb{R}$;
 \item \label{lc}$x(t)<0\; \forall\, t
\in \mathbb{R}$ if $\al(t)<0$.
\end{enumerate}
\end{lma}
\begin{proof}If $x(t)$ changes sign on $[0,T]$. then there is a $\tau$
such that $x(\tau)=0$. We may assume that $x'(\tau)\leq0 $.
Otherwise, if $x'(\tau)>0$, since $x(t)$ is a $T-$periodic function,
$x(t)$ has a successive zero $t_{0}$ on
 $ [\tau,\tau+T )$, such that  $x'(t_{0})\leq0$. In this case, let $\tau=t_{0}$. Without loss the
generality, we may assume that $\tau=0$ so that $x'(0)\leq0$. Let
$v(t)$ be initial value problem of the following equation

\begin{equation}
 \label{eqtf} y^{\prime \prime }-cy^{\prime  }+\al(t)y=0,
 \end{equation}

such that $v(0)=v'(0)-1=0$. by assumption that
$\al(t)\leq \frac{\pi^{2}}{%
T^{2}}+\frac{c^{2}}{4}=\lambda_{1}$ the first eigenvalue of
\[
y^{\prime \prime }-cy^{\prime  } +\lambda x=0, \, x(0)=x(T)=0.
\]
Therefore the equation
(\ref{eqtf}) disconjugate on $[0,T]$. Thus $v(t)>0 \forall t\in
(0,T]$. Multiplying (\ref{eqte}) by $v(t)$ subtracted from
(\ref{eqtf}) by $x(t)$, and integration by parts, we have that

\begin{equation}\label{eqtg}
v(T)x'(T)=vx'-xv'|_{0}^{T}=\int_{0}^{T}v(t)h(t)dt.\end{equation}

The left side of above equation is negative, while the right side is
positive, a contradiction.\end{proof}

\section{Proof of Theorem \ref{Theorema}}

Before giving the proof of Theorem \ref{Theorema}, let us list some results
that will be needed in the sequel.

Consider the following Duffing equation:
\begin{equation}
 \label{eqtha}x''+cx'+g(t,x)= h(t),
\end{equation}
  where $g(t,x)$
is
$T$-periodic in $t$ and differentiable in $x$.
\begin{lma}
\label{th1}Assume that
\begin{equation}
 g'_{x}(t,x)\ll
\left(\frac{2\pi}{T}\right)^{2}+\frac{c^2}{4}.\label{eqth1}
\end{equation}
\begin{enumerate}
\item\label{1a} The $T$-periodic solutions of \textup{(\ref{eqtha})} are
totally  ordered; \item \label{2b}\textup{(\ref{eqtha})} has a unique
$T$-periodic solution on
    $[u,v]$ if a $T$-periodic solution exists and $ g'_{x}(t,x)\gtrless 0$ on  $[u,v]$;
\item\label{3c} \textup{(\ref{eqtha})} cannot admit three distinct
$T$-periodic solutions in
    $[u,v]$ if $ g'_{x}(t,x)$ is  strictly increasing or strictly decreasing in
    $[u,v]$.
\end{enumerate}
\end{lma}

\begin{proof}The idea we use to prove the lemma is well known; we give the proof
here for completeness.

Let $x_{1}(t)$ and $x_{2}(t)$ be two distinct $T$-periodic solutions of (\ref{eqtha}),
and $u=x_{2}(t)-x_{1}(t)$ a nontrivial $T$-periodic solution of (\ref{eqtc}) with $\al(t)=
 \int_{0}^{1}g'(t, (1-s)x_{1}+sx_{2})\,ds$. Hence, the conclusion of (\ref{1a}) follows from
Lemma \ref{tlma}(\ref{l(1)}).
 Similarly, (\ref{2b}) follows from Lemma \ref{tlma}(\ref{l(3)}).
 Next, let $x_{1}(t)$, $x_{2}(t)$ and $x_{3}(t)$ be three distinct
 $T$-periodic solutions of (\ref{eqtha}) in the interval.
By (\ref{1a}), we may assume that $x_{1}(t)<x_{2}(t)<x_{3}(t)$.
 Setting $u_{i}=x_{i+1}(t)-x_{i}(t)$ for $i=1,2$, then $L_{\al_{i}}(u_{i})=0$ with
 $\al_{i}=[g(t,x_{i+1})-[g(t,x_{i})]/u_{i}$. The strict convexity of $g$ implies that
 $\al_{1}<\al_{2}$.
 By Lemma \ref{tlma}(\ref{l(2)}), we have $u_{1}(t)\equiv0$ or
 $u_{2}(t)\equiv0$, a contradiction in either case.
\end{proof}

 Now we turn to (\ref{aa}).
Let $Fx:=x''+cx'+ax-x^{3}-h(t)$ and let $B_{R}$ be a ball of radius $R$. We
have the following lemma.

\begin{lma}
\label{th2}For any positive number $a>0$, there is an $R>0$ large
enough that $\deg(F,B_{R},0)=-1$.
\end{lma}

\begin{proof} If $x(t)$ is a $T-$periodic solution of (\ref{aa}), let $t_{\max},\,t_{\min}$ be points
at which $x(t)$ achieves its maximum and minimum, respectively. Then
$x'(t_{\max})=0,\,x''(t_{\max})\leq0$, hence (\ref{aa})implies
\[ ax(t_{\max})-(x(t_{\max}))^{3}\geq h(t_{\max})>0
\]
so that \begin{equation} \label{eqthb} \max_{t}x(t)=x(t_{\max})\in
(-\infty,-\sqrt{a})\cup (0,\sqrt{a}).\end{equation}

 Similarly
\[ ax(t_{\min})-(x(t_{\min}))^{3}\leq h(t_{\min})<\parallel h\parallel _{\infty}
\]
so\[\min_{t}x(t)\geq C\] where $C$ is the negative root of
$g(C)=\parallel h\parallel _{\infty}$. Thus any $T-$periodic
solutions is bounded by \[C\leq x(t)\leq \sqrt{a},\] so the required
a priori bound is proved. The fact that the degree is $-1$ can also
be seen by using the homotopy \begin{equation} \label{eqthc}
x''+cx'+ax-x^{3}=\lambda h(t),\, 0\leq \lambda \leq 1.\end{equation}
For $\lambda =0$ is easy to compute directly that the degree is
$-1$, and since the priori bound above holds for all $\lambda \in
[0,1]$, we conclude that the degree is $-1$.
\end{proof}

Now we are ready to prove Theorem \ref{Theorema}.

First, we shall show that (\ref {aa}) has a unique  $T$-periodic
solution
 which is negative and unstable for $h(t)>h_{0} \;\forall\, t \in \mathbb{R}$.
It is obvious that $-R$ is a constant subsolution of (\ref {aa}) for
$R$ large enough. By the choice of $h_{0}$, it is easy to verify
that  $g(-2\sqrt{a/3})=h_{0}<h(t)$, so $b=-2\sqrt{a/3}$ is a
supersolution of (\ref {aa}).  By applying Lemma 3.2 in \cite{Zort},
there is a  $T$-periodic solution $x(t)$ of (\ref {aa}) such that
\[-R<x(t)<b.\]

Next, we  have to show that the solution obtained above is the only
$T$-periodic solution  of (\ref {aa}). Suppose there is another
$T$-periodic solution $y(t)$ of (\ref {aa}): then
$y(t)<b:=-2\sqrt{a/3}$. In fact, let let $t_{\max}$ be the point at
which $x(t)$ achieves its maximum, then \begin{equation}
\label{eqthcbis} ay(t_{\max})-(y(t_{\max}))^{3}\geq
h(t_{\max})>h_{0},
\end{equation} which forces that $ \max_{t}y(t)<b$, since $g(y)\leq h_{0}$ when $y>b$.

Now that both $x(t)$ and $y(t)$ are in $[-\infty,b]$, $g$ is
decreasing on this interval, so it follows from the second conclusion
of Lemma \ref{th1} that $x(t)\equiv y(t)$. Finally, since the
solution is unique, the local index of the $T$-periodic solution
is equal to the degree given by Lemma \ref{th2}, so the unique
$T$-periodic solution is unstable.

This completes the proof of the first part
of the theorem.

The idea of the proof of the second part is the same as that of
the first one, but here we have to estimate the solutions more
carefully.

 First, we note that (\ref {aa}) does not admit any sign-changing
 $T$-periodic solutions. This follows directly by rewriting
(\ref {aa}) as the following form:
\begin{equation}
 \label{eqthe}x''+cx'+q(t)x= h(t),
\end{equation}
  where
 $q(t)=g(x(t))/x(t)$ verifies the condition of Lemma \ref{lae}.
 Therefore, $x(t)$ is of constant sign.

 Next, we shall show the
 existence and uniqueness of a negative solution.

 The existence is
 evident, since $-\sqrt{a}$ and $b$ are constant sub- and supersolutions of
 (\ref{aa}) respectively in the case $0<h(t)<h_{0}$, so there is a
 $T$-periodic solution $x_{1}(t)$ of (\ref{aa}) such that
 $b<x_{1}(t)<-\sqrt{a}$.

 Letting $x(t)$ be any negative solution
 of (\ref {aa}), similar argument as the proof of Lemma \ref{th2}.
 from (\ref{eqthb}) we have \[\max_{t}x(t)<-\sqrt{a}.\] and from
 condition $0<h(t)<h_{0}$, we have

\[ ax(t_{\min})-(x(t_{\min}))^{3}\leq h(t_{\min})<h_{0}
\]
hence $\min x(t)>b$ since $g(x)\geq h_{0}$ for $x\leq b$.
  Therefore both $x_{1}(t)$ and $x(t)$
lie in the same interval $(b,-\sqrt{a})$. Moreover, $g$ is decreasing
on the interval, so Lemma \ref{th1} implies that $x_{1}(t)\equiv x(t)$.
This establishes the uniqueness of the negative solution and shows that $\operatorname
{ind}(x_{1}(t))=-1$.

Following the same reasoning as in the first part,  there is a
unique $T$-periodic solution $x_{3}(t)$ in $(\sqrt{a/3},\sqrt{a})$
with $\operatorname{ind}(x_{3}(t))=-1$. According to a formula  that
relates local index and topological degree, there is another
positive $T$-periodic solution $x_{2}(t)$. Since $g$ is concave on
$(0,\infty)$, by the third conclusion of Lemma \ref{th1} it is clear
that on $(0,\infty)$ (\ref{aa}) cannot admit more than two
solutions. We infer that the positive $T$-periodic solutions are
exactly two  and $\operatorname {ind}(x_{2}(t))=1$. This completes
the proof of the theorem. \md

We will finish the paper by proposing the following open questions
in the hope that the reader may solve it.\md

\begin{enumerate}
\item\label{ca} The main theorem characterizes the situation for positive forcing
 functions $h(t)$ whose graph does not cross the line $y=h_{0}$. The conjecture naturally arises
 that, for any forcing $h(t)$, or at least for any positive $h(t)$, there are at most three periodic
solutions (and generically either one or three solutions).\md
\item \label{cb}Is the condition \textup{(\ref{ab})}a sharp one for validity of the upper
bound on the number of solutions? In other words, assuming that
\[
a> \left(\frac{\pi}{T}\right)^{2}+\frac{c^{2}}{4},
\]
is it always possible to find  a forcing $h(t)$ with $h(t)>h_{0}$
such that (\ref{aa}) has more than one solution, and a forcing
$h(t)$ with $0<h(t)<h_{0}$ such that (\ref{aa}) has more than three
solutions?\md
\item\label{cc} We guess that
  the Duffing operator
$Fx:=x''+cx'+ax-x^{3}$ is globally equivalent to the cusp mapping in
the sense of Berger and Church \cite{Church} under the condition of
Theorem \ref{Theorema}. Up to now, we have not been able to give a
proof of it in the strict sense. So we state our question more
precisely below:

Let $\Sigma=\{\,x\in C_{T}^{2} \mid DF(x)v=0,\,v\neq0\,\}$ be the
singular set of $F$. Then $\Sigma$
 consists of the fold and cusp and $F(\Sigma)$ divides the space $C_{T}$ into two open components $A_{1}$, $A_{3}$.
Let $C$ be the subset that consists of
the cusp point. Then $F(C)$ is a co-dimension-one submanifold of
 $F(\Sigma)$ such that\md
\begin{enumerate}
\item\label{q1} (\ref {aa}) has a unique  $T$-periodic solution
 which is  unstable for  $h(t)\in A_{1}\cup F(C)$;
\item \label{q2}  (\ref{aa}) has  exactly three  ordered
$T$-periodic solutions if $h(t)\in A_{3}$.
 Moreover, the middle one is
asymptotically stable and the remaining two are unstable.
 \item\label{q3}(\ref{aa}) has  exactly two  ordered
$T$-periodic solutions for  $h(t)\in F(\Sigma)/F(C) $.
 Both of them are unstable.
\end{enumerate}\md
The main difficulty here is to verify that $F\colon\Sigma\rightarrow
F(\Sigma)$ is one to one.\end{enumerate}

\begin{acknowledgement}
The authors of this paper would like to thank the anonymous
referee for his/her careful reading of the paper and for many helpful
suggestions that improve the presentation of the ideas.
\end{acknowledgement}

\end{document}